\documentclass[12pt]{article}
\input isolatin1.sty
 \usepackage{a4}
 \usepackage{amsmath}
 \usepackage{amssymb}
\usepackage{graphicx}
\newtheorem{thm}{Theorem}[section]

\newtheorem{rem}[thm]{Remark}
\renewcommand{\a}{\alpha}

\title{Quantum Hilbert matrices and orthogonal polynomials}
\author{J{\o}rgen Ellegaard Andersen and Christian Berg}

\date{\today}
 \begin{document}
\maketitle

\begin{abstract} Using the notion of  quantum integers associated with
a complex number $q\neq 0$, we define the quantum Hilbert
matrix and various extensions. They are Hankel matrices
corresponding to certain little $q$-Jacobi polynomials when $|q|<1$, 
and for the special
value $q=(1-\sqrt{5})/(1+\sqrt{5})$ they are closely related to
Hankel matrices of reciprocal Fibonacci numbers called Filbert
matrices. We find a formula for the entries of the inverse quantum
Hilbert matrix.
\end{abstract}



\section{Introduction}

In \cite{Hi} Hilbert introduced the  {\it Hilbert matrices}
\begin{equation}\label{eq:hilbert}
\mathcal H_n=\left(1/(l+j+1)\right),\quad 0\le l,j\le n,\quad n=0,1,\ldots,
\end{equation}
 and found the following expression for their
determinants
\begin{equation}\label{hilbdet}
\det \mathcal H_n=\left(\prod_{k=1}^n (2k+1)
\tbinom{2k}{k}^2\right)^{-1},
\end{equation}
showing that they are the reci\-procal of integers. This
fact is also a consequence of the observation that the inverse matrices
$\mathcal H_n^{-1}$ have integer entries and Choi \cite{Ch} found the following
integer expression for them
\begin{equation}\label{eq:invhilb}
(\mathcal H_n^{-1})_{l,j}= 
(-1)^{l+j}(l+j+1)\tbinom{n+l+1}{n-j}\tbinom{n+j+1}{n-l}
\tbinom{l+j}{l}\tbinom{l+j}{j}.
\end{equation}

This formula was generalized in \cite[Theorem 4.1]{Be2} to the one
parameter extension $\mathcal H_n^{(\a)}=(\a/(l+j+\a))$ of the Hilbert matrices, where 
$\a>0$,
\begin{equation}\label{eq:invhilba}
(\mathcal H_n^{(\a)})^{-1}_{l,j}= 
(-1)^{l+j}\tfrac{l+j+\a}{\a}\tbinom{n+l+\a}{n-j}\tbinom{n+j+\a}{n-l}
\tbinom{l+j+\a-1}{l}\tbinom{l+j+\a-1}{j}.
\end{equation}
The idea of proof is to observe that the matrices are the Hankel matrices
of a moment problem. After having determined the corresponding orthogonal
polynomials and their kernel polynomials one uses the result that the
matrix of coefficients of the kernel polynomial is the inverse of the
Hankel matrix, see \cite[Theorem 2.1]{Be2}.
 
In \cite{Ri} Richardson noticed that the {\it Filbert
matrices}
\begin{equation}\label{eq:filbert}
\mathcal F_n=\left(1/F_{l+j+1}\right),\quad 0\le l,j\le n,\quad n=0,1,\ldots,
\end{equation}
where $F_n,n\ge 0$ is the sequence of Fibonacci numbers, have the
property that all elements of the inverse matrices are integers.
Richardson gave an explicit formula for the elements of the inverse
matrices and proved it using computer algebra. It is the special case
$\a=1,\sinh\theta=\tfrac12$ of (\ref{eq:filbinv}). The formula
shows a remarkable analogy with Choi's formula (\ref{eq:invhilb})
in the sense that one shall  replace the binomial coefficients
$\binom{n}{k}$ by the analogous {\it Fibonomial  coefficients}
\begin{equation}\label{eq:fibonomial}
\binom{n}{k}_{\mathbb F}=\prod_{j=1}^k\frac{F_{n-j+1}}{F_j},\quad 0\le
k\le n,
\end{equation}
with the usual convention that  empty products are defined as 1.
These coefficients are defined and studied in \cite{Kn} and are
integers.  We recall that the sequence
 of Fibonacci numbers is $F_0=0, F_1=1,\ldots,$ with the recursion formula
$F_{n+1}=F_n+F_{n-1},\;n\ge 1$.

The Hilbert matrices are the Hankel matrices $(s_{l+j})$ corresponding to the
moment sequence 
\begin{equation}\label{eq:momlegendre}
s_n=1/(n+1)=\int_0^1 x^n\,dx,\quad n\ge 0,
\end{equation}
and the corresponding orthogonal polynomials are the Legendre
polynomials for the interval $[0,1]$, see \cite[Section 7.7]{A:A:R}.
This was used in \cite{Be2} and in the survey paper \cite{Be} to prove
Choi's formula.

In \cite{Be2} it was also proved that for $\a\in\mathbb N$ the
sequence $(F_\a/F_{n+\a})_{n\ge 0}$ is the
moment sequence of a signed measure of total mass one and 
the corresponding orthogonal polynomials were identified as special little
$q$-Jacobi polynomials corresponding to the value
$q=(1-\sqrt{5})/(1+\sqrt{5})$. From this Richardson's formula
(\ref{eq:filbinv}) for the elements of the inverse Filbert matrices
was derived as the case $\a=1$ of \cite[(29)]{Be2}.

The results about Fibonacci numbers have been extended by Ismail in 
\cite{Is1} to a one-parameter family of sequences
$(F_n(\theta))_{n\ge 0}, \theta>0$,
determined by the recursion
\begin{equation}\label{eq:genfib}
F_{n+1}(\theta)=2\sinh\theta F_n(\theta)+F_{n-1}(\theta),\; n\ge 1,\quad
F_0(\theta)=0,F_1(\theta)=1.
\end{equation}
When $\sinh\theta=\tfrac12$ we have $F_n(\theta)=F_n$, and when
$2\sinh\theta$ is a positive integer then all $F_n(\theta)$ are
integers.

The purpose of this paper is to give a common generalization of all
these results by the use of {\it quantum integers} defined by
\begin{equation}\label{eq:quantum}
[n]_q=\frac{q^{n/2}-q^{-n/2}}{q^{1/2}-q^{-1/2}},\quad n=0,1,\ldots
\end{equation}
where $q\in \mathbb C\setminus\{0\}$.
To make the definition precise we consider $q\to [n]_q$ as a
holomorphic function in the cut plane $\mathbb C\setminus]-\infty,0]$
extended to $q\in \left]-\infty,0\right[$ by
$$
[n]_q=\lim_{\varepsilon\to 0^+} [n]_{q+i\varepsilon}.
$$
Clearly $[n]_q=n$ for $n=0,1$ and $[n]_1=n$ for all $n\in\mathbb N$.
The quantum integer $[n]_q$ vanishes, when $q$ is an $n$'th root of
unity different from 1.

Assume now $0<|q|<1$.  We  consider the complex measures in the
 complex plane depending on
$\a\in\mathbb N$ and defined by
\begin{equation}\label{eq:compmea}
\mu^{(\a)}(q)=(1-q^\a)\sum_{k=0}^\infty q^{k\a}\delta_{q^{k+1/2}},
\end{equation}
where $\delta_c$ is the Dirac measure concentrated at $c\in\mathbb
C$. It is easy to see the following $q$-analogue of (\ref{eq:momlegendre})
\begin{equation}\label{eq:moments}
\frac{1}{[n+1]_q}=\int x^n\,d\mu^{(1)}(q)(x),\quad n\ge 0,
\end{equation}
and in general
\begin{equation}\label{eq:qamoments}
\frac{[\a]_q}{[n+\a]_q}=\int x^n\,d\mu^{(\a)}(q)(x),\quad n\ge 0.
\end{equation}
The Hankel matrices corresponding to $\mu^{(1)}(q)$ are defined as
\begin{equation}\label{eq:qhankel}
\mathcal H_n(q)=\left(\frac{1}{[l+j+1]_q}\right),\quad 0\le l,j\le n,\quad n=0,1,\ldots,
\end{equation}
and are called the {\it quantum Hilbert matrices}, but we will also
consider the {\it generalized quantum Hilbert matrices} ($\a\in\mathbb N$)
\begin{equation}\label{eq:qahankel}
\mathcal H_n^{(\a)}(q)=\left(\frac{[\a]_q}{[l+j+\a]_q}\right),\quad 0\le l,j\le n,
\quad n=0,1,\ldots.
\end{equation}
These matrices are well-defined  for non-zero complex numbers $q$
which are not roots of unity of order $\le 2n+\a$.
When $0<q<1$ the measure $\mu^{(1)}(q)$ is a probability measure on $[0,1]$
which converges weakly to the Lebesgue measure on $[0,1]$ for $q\to 1$.
The quantum Hilbert matrices $\mathcal H_n(q)$ converge to the ordinary Hilbert
matrices when $q\to 1$.

We prove in section 2 that the generalized quantum Hilbert matrices are
regular and find a formula for the elements of the inverse matrix, see
Theorem \ref{thm:inverseqhilb}.
This is a $q$-analogue of Choi's formula when $\a=1$.
The proof uses the method of \cite{Be2} by identifying the
orthogonal polynomials of $\mu^{(\a)}(q)$ as little $q$-Jacobi polynomials.

In section 3 we consider the special values
$q=-e^{-2\theta},\;\theta>0$, which for $\sinh\theta=\tfrac12$ gives 
$q=(1-\sqrt{5})/(1+\sqrt{5})$.
We prove that $\mathcal H_n^{(\a)}(q)$ is unitarily related to
 the {\it generalized Filbert matrix}
\begin{equation}\label{eq:genfilbert}
\mathcal F_n^{(\a)}(\theta)=(F_{\a}(\theta)/F_{l+j+\a}(\theta)),
\quad 0\le l,j\le n,\quad n=0,1,\ldots
\end{equation}
in the sense that
\begin{equation}\label{eq:quantunifilbert}
\mathcal H_n^{(\a)}(q)=U_n\mathcal F_n^{(\a)}(\theta) U_n,
\end{equation}
where $U_n$ is a unitary diagonal matrix with diagonal elements
$i^l,l=0,1,\ldots,n$.
This makes it possible to deduce Richardson's formula 
for the elements of the inverse  of $\mathcal F_n$ and its generalizations
$\mathcal F_n^{(\a)}(\theta)$ from the $q$-analogue of 
Choi's formula.
 
\section{Quantum Hilbert matrices}

The little $q$-Jacobi
polynomials are given by
\begin{equation}\label{eq:lqJacobi}
p_n(x;a,b;q)={}_2\phi_1\left(\begin{matrix}q^{-n},abq^{n+1}\\aq\end{matrix};q,xq\right),
\end{equation}
where $a,b$ are complex parameters satisfying
$|a|,|b|\le 1$.

In \cite[Section 7.3]{G:R} one finds a discussion of the little
$q$-Jacobi polynomials defined in (\ref{eq:lqJacobi}), and it is
proved that
\begin{equation}\label{eq:ortJacobi}
\sum_{k=0}^\infty p_n(q^k;a,b;q)p_m(q^k;a,b;q)\frac{(bq;q)_k}{(q;q)_k}
(aq)^k=\frac{\delta_{n,m}}{h_n(a,b;q)},
\end{equation}
where
\begin{equation}\label{eq:norm}
h_n(a,b;q)=\frac{(abq;q)_n(1-abq^{2n+1})(aq;q)_n(aq;q)_\infty}
{(q;q)_n(1-abq)(bq;q)_n(abq^2;q)_\infty}(aq)^{-n}.
\end{equation}
In \cite{G:R} it is assumed that $0<q,aq<1$, but the derivation shows
that it holds for $0<|q|<1,|a|\le 1, |b|\le 1$, in particular in the
case of interest here: $a=q^{\a-1},\a\in\mathbb N,b=1$, in the case of which we get  
\begin{equation}\label{eq:ortJacobispec}
\sum_{k=0}^\infty p_n(q^k;q^{\a-1},1;q)p_m(q^k;q^{\a-1},1;q)
q^{\a k}=\delta_{n,m}\frac{q^{\a n}(q;q)_n^2}{(q^{\a};q)_n^2(1-q^{2n+\a})}.
\end{equation}

The Gaussian $q$-binomial coefficients
$$
\left[\begin{matrix}n\\k\end{matrix}\right]_q=
\frac{(q;q)_n}{(q;q)_k(q;q)_{n-k}}
$$
are polynomials in $q$. Using the quantum integers  (\ref{eq:quantum})
we get
$$
\left[\begin{matrix}n\\k\end{matrix}\right]_q=\prod_{j=1}^k
\frac{[n-j+1]_q q^{(n-j)/2}}{[j]_q q^{(j-1)/2}}=
q^{k(n-k)/2} \prod_{j=1}^k \frac{[n-j+1]_q}{[j]_q},
$$
and introducing the notation 
$$
[n]_q!=\prod_{k=1}^n [k]_q,\quad [0]_q!=1,
$$
we define the {\it quantum binomial coefficients} by
\begin{equation}\label{eq:quantumbin}
\binom{n}{k}_q=
 \prod_{j=1}^k \frac{[n-j+1]_q}{[j]_q}=\frac{[n]_q!}{[k]_q![n-k]_q!},
\end{equation}
hence
\begin{equation}\label{eq:gaussquantum}
\left[\begin{matrix}n\\k\end{matrix}\right]_q=q^{k(n-k)/2}\binom{n}{k}_q.
\end{equation}

This formula shows that $\binom{n}{k}_q$ is a holomorphic function of
$q$ in the cut plane $\mathbb C\setminus ]-\infty,0]$ with a
continuous extension to the upper part of the cut.

Defining
\begin{equation}\label{eq:fibpol}
p_n^{(\a)}(q;x):=\binom{n+\a-1}{n}_q p_n(x/q^{1/2};q^{\a-1},1;q),
\end{equation}
some calculation leads to the simple expression
\begin{equation}\label{eq:fibpol2}
p_n^{(\a)}(q;x)=\sum_{j=0}^n \binom{n}{j}_q
\binom{n+j+\a-1}{n}_q(-1)^j x^j,
\end{equation}
which is a $q$-analogue of the polynomials $r_n^{(\a)}$ of \cite{Be2}.

The equation (\ref{eq:ortJacobispec}) can be written
\begin{equation}\label{eq:ortfinal}
\int p_n^{(\a)}(q;x)p_m^{(\a)}(q;x)\,d\mu^{(\a)}(q)(x)=\delta_{n,m}
\frac{[\a]_q}{[2n+\a]_q}
\end{equation}
showing that
$$
P_n^{(\a)}(q;x)=([2n+\a]_q/[\a]_q)^{1/2}p_n^{(\a)}(q;x),
$$
are orthonormal polynomials. The corresponding kernel polynomials are
defined by
\begin{equation}\label{eq:qkernel}
K_n^{(\a)}(q;x,y)=\sum_{k=0}^n P_k^{(\a)}(q;x)P_k^{(\a)}(q;y)=
\sum_{k=0}^n\left([2k+\a]_q/[\a]_q\right)p_k^{(\a)}(q;x)p_k^{(\a)}(q;y).
\end{equation}
While $P_n^{(\a)}(q;x)$ depends on the choice of a square root, the kernel
polynomials $K_n^{(\a)}(q;x,y)$ are independent of this choice. Since the 
orthogonal polynomials exist with respect to $\mu^{(\a)}(q)$ and the
moments are given by (\ref{eq:qamoments}), it follows by
\cite[Theorem 3.1]{Chi} that $\mathcal H_n^{(\a)}(q)$ is regular for each $n$.
For $q\to 1$ and $\a=1$ the expressions (\ref{eq:fibpol2})--(\ref{eq:qkernel})
tend to classical formulas for Legendre polynomials for $[0,1]$.

Writing
$$
K_n^{(\a)}(q;x,y)=\sum_{l,j=0}^n a^{(n)}_{l,j}(q;\a) x^ly^j,
$$
it follows by (\ref{eq:fibpol2}) that the coefficients $a^{(n)}_{l,j}(q;\a)$
are given by
\begin{equation}\label{eq:nyformel} 
 a^{(n)}_{l,j}(q;\a)=(-1)^{l+j}\sum_{k=\max{(l,j)}}^n
\tfrac{[2k+\a]_q}{[\a]_q}\tbinom{k}{l}_q\tbinom{k}{j}_q\tbinom{k+l+\a-1}{k}_q
\tbinom{k+j+\a-1}{k}_q.
\end{equation} 
\begin{thm}\label{thm:inverseqhilb} The $l,j$'th element of the inverse
  matrix of the generalized quantum Hilbert matrix $\mathcal
  H_n^{(\a)}(q)$ defined in (\ref{eq:qahankel}) is given as
\begin{equation}\label{eq:qchoi} 
 (-1)^{l+j}\tfrac{[l+j+\a]_q}{[\a]_q}\tbinom{n+l+\a}{n-j}_q
\tbinom{n+j+\a}{n-l}_q\tbinom{l+j+\a-1}{l}_q\tbinom{l+j+\a-1}{j}_q.
\end{equation}
 Furthermore,
\begin{equation}\label{eq:qhilb}
\det \mathcal H_n^{(\a)}(q)=[\a]_q^n\left(\prod_{k=1}^n [2k+\a]_q
\tbinom{2k+\a-1}{k}_q^2\right)^{-1}.
\end{equation}
\end{thm}

{\it Proof}. It is a general fact that the coefficients $a^{(n)}_{l,j}(q)$
of the kernel polynomial are the entries of the inverse of the Hankel
matrix, cf. \cite[Theorem 2.1]{Be2}.
 
 Let $R(n;l,j)$ denote the number given in the right-hand side of
(\ref{eq:qchoi}), and define 
$$
C(k;l,j)=(-1)^{l+j}\tfrac{[2k+\a]_q}{[\a]_q}\tbinom{k}{l}_q\tbinom{k}{j}_q
\tbinom{k+l+\a-1}{k}_q\tbinom{k+j+\a-1}{k}_q, \quad k\ge l,j.
$$
We shall prove that
\begin{equation}\label{eq:induct}
R(n;l,j)=\sum_{k=\max(l,j)}^n C(k;l,j)
\end{equation}
by induction in $n$, and can assume $l\ge j$ without loss of
generality.
 The equation (\ref{eq:induct}) is easy for $n=k=l$
and is left to the reader. We shall establish the induction step
\begin{equation}\label{eq:induc}
R(n+1;l,j)-R(n;l,j)=C(n+1;l,j).
\end{equation}
The left-hand side of this expression can  be written
$$
(-1)^{l+j}\tfrac{[l+j+\a]_q}{[\a]_q}\tbinom{l+j+\a-1}{l}_q
\tbinom{l+j+\a-1}{j}_qT,
$$
where
$$
T=\tbinom{n+l+\a+1}{n+1-j}_q\tbinom{n+j+\a+1}{n+1-l}_q-\tbinom{n+l+\a}{n-j}_q
\tbinom{n+j+\a}{n-l}_q
$$
$$
=\tfrac{([n+l+\a]_q\cdots[l+j+\a+1]_q)([n+j+\a]_q\cdots[l+j+\a+1]_q)}
{[n+1-j]_q![n+1-l]_q!}\,\cdot
$$
$$
\left\{[n+l+\a+1]_q[n+j+\a+1]_q-[n+1-j]_q[n+1-l]_q\right\}.
$$
The quantity in braces equals $[2n+2+\a]_q[l+j+\a]_q$, and now it is
easy to complete the proof of (\ref{eq:induc}).

From the general theory of orthogonal polynomials,
cf. \cite{Ak},\cite{Chi},\cite{Is}, it is  known that the leading coefficient of the
orthonormal polynomial $P_n^{(\a)}(q;x)$ is $\sqrt{D_{n-1}/D_n}$, where
$$
D_n=\det\mathcal H_n^{(\a)}(q).
$$
From (\ref{eq:fibpol2}) and (\ref{eq:ortfinal}) we then get
$$
D_{n-1}/D_n=([2n+\a]_q/[\a]_q)
\tbinom{2n+\a-1}{n}_q^2,
$$
hence
$$
\frac{1}{D_n}=\prod_{k=1}^n
\frac{D_{k-1}}{D_{k}}=
\prod_{k=1}^n ([2k+\a]_q/[\a]_q) \tbinom{2k+\a-1}{k}_q^2,
$$
and (\ref{eq:qhilb}) follows.
 \quad$\square$

\begin{rem} {\rm By analytic continuation, the formulas
(\ref{eq:qchoi}) and (\ref{eq:qhilb}) of Theorem \ref{thm:inverseqhilb}
are valid for $q\ne 0$ which is not a root of unity of order $\le 2n+\a$.}  
\end{rem}

\section{Filbert matrices}

In this section we specialize to $q=-e^{-2\theta}$ for $\theta>0$. It
is easy to see that the unique solution $F_n(\theta),n\ge 0$ to 
(\ref{eq:genfib}) is given by
\begin{equation}\label{eq:genfib1}
F_n(\theta)=\frac{e^{n\theta}-(-1)^ne^{-n\theta}}{e^\theta+e^{-\theta}},
\end{equation}
which is Ismail's definition in \cite{Is1}, hence
\begin{equation}\label{eq:quantumfilbert}
[n]_q=(-i)^{n-1}F_n(\theta).
\end{equation}
For $\theta=\theta_0>0$ such that $\sinh\theta_0=\tfrac12$ we get
 $q=(1-\sqrt{5})/(1+\sqrt{5})$ and $F_n=F_n(\theta_0)$. 
For information about Fibonacci numbers, see
\cite{Kn},\cite{Ko}.
Ismail \cite{Is1} also considered the generalized Fibonomial
coefficients
\begin{equation}\label{eq:genfibonomial}
\binom{n}{k}_{\mathbb F(\theta)}=\prod_{j=1}^k\frac{F_{n-j+1}(\theta)}
{F_j(\theta)},\quad 0\le
k\le n,
\end{equation}
with the usual convention that empty products are 1, and gave the
following recursion
\begin{equation}\label{eq:genfibrec}
\binom{n}{k}_{\mathbb F(\theta)}=F_{k-1}(\theta)
\binom{n-1}{k}_{\mathbb F(\theta)}+F_{n-k+1}(\theta)
\binom{n-1}{k-1}_{\mathbb F(\theta)},\quad n>k\ge 1,
\end{equation}
which shows that they are integers when $2\sinh\theta$ is an integer.
Using (\ref{eq:quantumfilbert}) we next get that the quantum binomial
 coefficients can be expressed by
the Fibonomial coefficients of (\ref{eq:genfibonomial}) as
\begin{equation}\label{eq:qbinom}
\binom{n}{k}_q=(-i)^{k(n-k)}\binom{n}{k}_{\mathbb F(\theta)},
\end{equation}
hence for $0\le l,j\le n, \a\in\mathbb N$
\begin{equation}\label{eq:qbinom1}
\frac{[\a]_q}{[l+j+\a]_q}=i^{l+j}\frac{F_\a(\theta)}{F_{l+j+\a}(\theta)}.
\end{equation}
Letting $U_n$ denote the unitary diagonal $(n+1)\times
(n+1)$-matrix with $l$'th
diagonal element equal to $i^l, l=0,1,\ldots,n$, then formula
 (\ref{eq:qbinom1}) implies
\begin{equation}\label{matrix}
\mathcal H_n^{(\a)}(q)=U_n\mathcal F_n^{(\a)}(\theta) U_n,\quad 
\mathcal F_n^{(\a)}(\theta)^{-1}=U_n
\mathcal H_n^{(\a)}(q)^{-1}U_n.
\end{equation} 
This leads to a new proof of Berg's and Ismail's generalizations of 
Richardson's formula, cf. \cite{Be2}, \cite{Is1}.

\begin{thm}\label{thm:genfilbinv} Let $A$ be the matrix 
$(1/F_{l+j+\a}(\theta)),0\le l,j\le n$. Then $A^{-1}$ has the entries
\begin{equation}\label{eq:filbinv}
(-1)^{n(l+j+\a)-\binom{l}{2}-\binom{j}{2}} F_{l+j+\a}(\theta)
\tbinom{n+l+\a}{n-j}_{\mathbb F(\theta)}
\tbinom{n+j+\a}{n-l}_{\mathbb F(\theta)}
\tbinom{l+j+\a-1}{l}_{\mathbb F(\theta)}
\tbinom{l+j+\a-1}{j}_{\mathbb F(\theta)},
\end{equation}
and
\begin{equation}\label{eq:det}
\det A=(-1)^{\a\binom{n+1}{2}}\left(F_\a(\theta)\prod_{k=1}^n F_{2k+\a}(\theta)
\tbinom{2k+\a-1}{k}^2_{\mathbb F(\theta)}\right)^{-1}.
\end{equation} 
\end{thm}

{\it Proof.} We use 
$$
\det U_n=i^{\binom{n+1}{2}},\quad \det\mathcal
H_n^{(\a)}(q)=(-1)^{\binom{n+1}{2}}
\det\mathcal F_n^{(\a)}(\theta)
$$
and the formulas (\ref{eq:quantumfilbert}) and (\ref{eq:qbinom}) to
make the calculation, the only non-obvious thing being the sign in the
two formulas. In the proof of (\ref{eq:filbinv}) we get the following
sign for the $lj$'th element of $A^{-1}$:
\begin{eqnarray*}
\lefteqn{
i^{l+j}(-1)^{l+j}(-i)^{l+j+(n-j)(l+j+\a)+(n-l)(l+j+\a)+l(j+\a-1)+j(l+\a-1)}}\\
&=&(-1)^{l+j}(-i)^{2n(l+j+\a)-l(l+1)-j(j+1)}=(-1)^{n(l+j+\a)-\binom{l}{2}-\binom{j}{2}}.
\end{eqnarray*}
In the proof of (\ref{eq:det}) we get the sign
\begin{eqnarray*}
\lefteqn{(-1)^{\binom{n+1}{2}}(-i)^{(\a-1)n}\left(\prod_{k=1}^n
  (-i)^{2k+\a-1+2k(k+\a-1)}\right)^{-1}}\\
&=&
(-1)^{\binom{n+1}{2}}\left(\prod_{k=1}^n
  (-i)^{2k(k+1)+2k(\a-1)}\right)^{-1}
=(-1)^{\binom{n+1}{2}}\left(\prod_{k=1}^n (-1)^{k(\a-1)}\right)^{-1}\\
&=&(-1)^{\a\binom{n+1}{2}}.
\end{eqnarray*}

 $\quad\square$

\medskip
J{\o}rgen Ellegaard Andersen, Department of Mathematics, University of
Aarhus, Ny Munkegade,
DK-8000 Aarhus C, Denmark. email: andersen@imf.au.dk

\medskip
Christian Berg, Department of Mathematics, University of Copenhagen, 
Universitetsparken 5, DK
2100 Copenhagen {\O}, Denmark. email: berg@math.ku.dk 
 \end{document}